\documentclass[english,12pt]{smfart}
\usepackage{amssymb}
\usepackage{amsfonts}
\usepackage{amscd}
\usepackage[T1]{fontenc}
\usepackage[all]{xypic}

\swapnumbers

\def\limproj{\mathop{\oalign{lim\cr
\hidewidth$\longleftarrow$\hidewidth\cr}}}

\def\kvar{{K_0 ({\rm Var}_k)}}

\def\Gr{{\rm Gr}}

\let\cal\mathcal

\def\AA{{\mathbf A}}

\def\CC{{\mathbf C}}
\def\DD{{\mathbf D}}

\def\FF{{\mathbf F}}

\def\LL{{\mathbf L}}

\def\NN{{\mathbf N}}

\def\QQ{{\mathbf Q}}

\def\ZZ{{\mathbf Z}}

\def\cI{{\mathcal I}}

\def\cL{{\mathcal L}}

\def\cN{{\mathcal N}}
\def\cO{{\mathcal O}}

\def\cR{{\mathcal R}}
\def\cS{{\mathcal S}}

\def\cU{{\mathcal U}}

\def\cX{{\mathcal X}}
\def\cY{{\mathcal Y}}

\mathchardef\alphag="7C0B
\mathchardef\betag="7C0C
\mathchardef\gammag="7C0D
\mathchardef\deltag="7C0E
\mathchardef\varepsilong="7C22
\mathchardef\varphig="7C27
\mathchardef\psig="7C20
\mathchardef\zetag="7C10
\mathchardef\epsilong="7C0F
\mathchardef\rhog="7C1A
\mathchardef\taug="7C1C
\mathchardef\upsilong="7C1D
\mathchardef\iotag="7C13
\mathchardef\thetag="7C12
\mathchardef\pig="7C19
\mathchardef\sigmag="7C1B
\mathchardef\etag="7C11
\mathchardef\omegag="7C21
\mathchardef\kappag="7C14
\mathchardef\lambdag="7C15
\mathchardef\mug="7C16
\mathchardef\xig="7C18
\mathchardef\chig="7C1F
\mathchardef\nug="7C17
\mathchardef\varthetag="7C23
\mathchardef\varpig="7C24
\mathchardef\varrhog="7C25
\mathchardef\varsigmag="7C26
\mathchardef\Omegag="7C0A
\mathchardef\Thetag="7C02
\mathchardef\Sigmag="7C06
\mathchardef\Deltag="7C01
\mathchardef\Phig="7C08
\mathchardef\Gammag="7C00
\mathchardef\Psig="7C09
\mathchardef\Lambdag="7C03
\mathchardef\Xig="7C04
\mathchardef\Pig="7C05
\mathchardef\Upsilong="7C07

\newtheorem{theorem}[subsubsection]{Theorem}
\newtheorem{lem}[subsubsection]{Lemma}
\newtheorem{cor}[subsubsection]{Corollary}
\newtheorem{prop}[subsubsection]{Proposition}
\newtheorem{problem}[subsubsection]{Problem}

\theoremstyle{definition}
\newtheorem{definition}[subsubsection]{Definition}

\newtheorem{def-prop}[subsubsection]{Proposition-Definition}
\newtheorem{def-theorem}[subsubsection]{Theorem-Definition}

\theoremstyle{remark}
\newtheorem{remark}[subsubsection]{Remark}

\theoremstyle{plain}

\numberwithin{equation}{subsection}

\def\boxit#1#2{\setbox1=\hbox{\kern#1{#2}\kern#1}%
\dimen1=\ht1 \advance\dimen1 by #1
\dimen2=\dp1 \advance\dimen2 by #1
\setbox1=\hbox{\vrule height\dimen1 depth\dimen2\box1\vrule}%
\setbox1=\vbox{\hrule\box1\hrule}%
\advance\dimen1 by .4pt \ht1=\dimen1
\advance\dimen2 by .4pt \dp1=\dimen2 \box1\relax}

\renewcommand{\theequation}{\thesubsection.\arabic{equation}}

\let\cal\mathcal

\def\AA{{\mathbf A}}

\def\CC{{\mathbf C}}
\def\DD{{\mathbf D}}

\def\FF{{\mathbf F}}

\def\LL{{\mathbf L}}

\def\NN{{\mathbf N}}

\def\QQ{{\mathbf Q}}

\def\ZZ{{\mathbf Z}}

\def\cI{{\mathcal I}}

\def\cL{{\mathcal L}}

\def\cN{{\mathcal N}}
\def\cO{{\mathcal O}}

\def\cR{{\mathcal R}}
\def\cS{{\mathcal S}}

\def\cU{{\mathcal U}}

\def\cX{{\mathcal X}}
\def\cY{{\mathcal Y}}

\mathchardef\alphag="7C0B
\mathchardef\betag="7C0C
\mathchardef\gammag="7C0D
\mathchardef\deltag="7C0E
\mathchardef\varepsilong="7C22
\mathchardef\varphig="7C27
\mathchardef\psig="7C20
\mathchardef\zetag="7C10
\mathchardef\epsilong="7C0F
\mathchardef\rhog="7C1A
\mathchardef\taug="7C1C
\mathchardef\upsilong="7C1D
\mathchardef\iotag="7C13
\mathchardef\thetag="7C12
\mathchardef\pig="7C19
\mathchardef\sigmag="7C1B
\mathchardef\etag="7C11
\mathchardef\omegag="7C21
\mathchardef\kappag="7C14
\mathchardef\lambdag="7C15
\mathchardef\mug="7C16
\mathchardef\xig="7C18
\mathchardef\chig="7C1F
\mathchardef\nug="7C17
\mathchardef\varthetag="7C23
\mathchardef\varpig="7C24
\mathchardef\varrhog="7C25
\mathchardef\varsigmag="7C26
\mathchardef\Omegag="7C0A
\mathchardef\Thetag="7C02
\mathchardef\Sigmag="7C06
\mathchardef\Deltag="7C01
\mathchardef\Phig="7C08
\mathchardef\Gammag="7C00
\mathchardef\Psig="7C09
\mathchardef\Lambdag="7C03
\mathchardef\Xig="7C04
\mathchardef\Pig="7C05
\mathchardef\Upsilong="7C07

\DeclareMathOperator*{\Spec}{Spec}
\DeclareMathOperator*{\Specf}{Spf}

\def\ord{{\rm ord}}
\def\Jac{{\rm Jac}}

\begin{document}

\title[Motivic integration 
on smooth rigid varieties]{Motivic integration
on smooth rigid varieties and invariants of degenerations}

\author{Fran\c cois Loeser}

\address{{\'E}cole Normale Sup{\'e}rieure,
D{\'e}partement de math{\'e}matiques et applications,
45 rue d'Ulm,
75230 Paris Cedex 05, France
(UMR 8553 du CNRS)}
\email{Francois.Loeser@ens.fr}
\urladdr{http://www.dma.ens.fr/~loeser/}

\author{Julien Sebag}

\address{{\'E}cole Normale Sup{\'e}rieure,
D{\'e}partement de math{\'e}matiques et applications,
45 rue d'Ulm,
75230 Paris Cedex 05, France
(UMR 8553 du CNRS)}
\email{Julien.Sebag@ens.fr}


\maketitle

\section{Introduction}In the last years, motivic integration
has shown to be a quite powerful tool in producing new invariants 
in birational geometry of algebraic varieties
over a field $k$, say  of characteristic zero, cf. \cite{K}\cite{Bat}\cite{Ba3}\cite{arcs}\cite{MK}\cite{barc}.
Let us explain the basic idea behind such results.
If $h : Y \rightarrow X$ is a proper birational morphism
between $k$-algebraic varieties, the induced morphism
$\cL (Y) \rightarrow \cL (X)$ between arc spaces (cf. \cite{arcs})
is an isomorphism outside subsets of infinite codimension.
By a fundamental 
change of variable formula,  
motivic integrals on $\cL (X)$ may be computed on $\cL (Y)$ 
when $Y$ is smooth.

In the present paper we develop similar ideas in the somewhat dual
situation of degenerating families over
complete discrete valuation 
rings with perfect residue field, for which
rigid geometry appears to be a natural framework.
More
precisely, let $R$ be a complete discrete valuation 
ring with fraction field $K$ and perfect
residue field $k$.
We construct a theory of motivic integration
for smooth\footnote{The extension to singular rigid spaces when $K$ is of
characteristic zero will be considered in a separate publication.}
rigid $K$-spaces, always assumed to be
quasi-compact and separated.
Let $X$ be a smooth rigid $K$-space of dimension $d$.
Our construction 
assigns 
to a gauge form $\omega$ on $X$, {\it i.e.} a 
nowhere vanishing differential form of degree $d$
on $X$, 
an integral $\int_{X} \omega d \tilde \mu$
with value
in the ring 
$\kvar_{\rm loc}$. Here $\kvar_{\rm loc}$ is the localization
with respect to the class of the affine line of the Grothendieck group 
of algebraic varieties
over $k$. In concrete terms, two varieties over $k$ define
the same class in $\kvar_{\rm loc}$ if they become isomorphic
after cutting them into locally closed pieces
and stabilization by product with a power of the affine line.
More generally, if $\omega$ is a differential form of degree $d$
on $X$, we define an integral 
$\int_{X} \omega d  \mu$ with value in the ring 
$\widehat \kvar$, 
which is the completion of 
$\kvar_{\rm loc}$ with respect to the filtration by 
virtual dimension (see \S\kern .15em \ref{GK}).
The construction is done by viewing $X$ as the generic fibre
of some formal $R$-scheme $\cX$. 
To such a 
formal $R$-scheme, by mean of 
the Greenberg functor $\cX \mapsto \Gr (\cX)$ one
associates a certain $k$-scheme $\Gr (\cX)$, which,
when $R = k [[t]]$,
and $\cX$ is the formal  completion of $X_0 \otimes k [[t]]$, for $X_0$ an 
algebraic variety over $k$, is nothing else that
the arc space
$\cL (X_0)$ considered before. We may then use the 
general theory of motivic integration
on schemes
$\Gr (\cX)$ which is developed in \cite{sebag}. 
Of course, for the construction
to work one needs to check that it is independent of the chosen model.
This is done by using two main ingredients: 
the theory of weak N{\'e}ron models developed in
\cite{neron} and \cite{bs}, and the analogue for schemes of
the form
$\Gr (\cX)$
of the
change of variable formula  which is proven
in
\cite{sebag}. In fact the theory
of weak N{\'e}ron models really 
pervades the whole paper and some parts of the book \cite{neron}
were crying out for their use in motivic integration\footnote{cf.
the remark at the bottom of p.105 in \cite{neron}.}.

As an application 
of our theory, we are able to assign in a canonical way
to any smooth quasi-compact and separated
rigid  $K$-space
$X$  an element $\lambda (X)$ in
the quotient ring $\kvar_{\rm loc} /(\LL - 1) \kvar_{\rm loc}$,
where $\LL$ stands for the class of the affine line.  When $X$ admits a formal
$R$-model with good reduction, $\lambda (X)$ is just the class of the
fibre of that model. More generally, if $\cU$ is a weak N{\'e}ron model of $X$,
$\lambda (X)$ is equal to the class of the special fibre of $\cU$ in
$\kvar_{\rm loc} /(\LL - 1) \kvar_{\rm loc}$. In particular it follows
that this class is independent of the choice
of the 
weak N{\'e}ron model $\cU$.

This result can be viewed as a rigid analogue of a result of Serre
concerning compact smooth 
locally analytic varieties over a local field. To such a variety $\tilde X$,
Serre associates in \cite{serre}, using classical
$p$-adic integration,
an invariant $s (\tilde X)$ in the
ring $\ZZ /(q - 1) \ZZ$, where $q$ denotes the
cardinality of the finite field
$k$. Counting rational points in $k$ yields a canonical morphism
$\kvar_{\rm loc} /(\LL - 1) \kvar_{\rm loc} \rightarrow \ZZ /(q - 1) \ZZ$,
and we show that the image by this morphism
of our motivic invariant $\lambda (X)$ of a smooth rigid $K$-space
$X$
is equal to the Serre invariant
of the underlying 
locally analytic
variety.

Unless making additional assumptions on $X$,
one cannot hope to lift our
invariant $\lambda (X)$ to a class in the Grothendieck ring
$\kvar_{\rm loc}$, which would be a substitute for the class
of the special fibre of {\bf the} N{\'e}ron model, when such a 
N{\'e}ron model happens
to exist. In the particular situation
where $X$ is the analytification of a 
Calabi-Yau variety over $K$, \textit{i.e.}
a smooth projective algebraic variety over $K$ of pure dimension $d$
with $\Omega^d_X$ trivial, this can be achieved: one can attach to
$X$ a canonical element of $\kvar_{\rm loc}$, which,
if $X$ admits a proper and smooth $R$-model $\cX$,
is equal to the class
of the special fibre $\cX_0$ in $\kvar_{\rm loc}$. In particular,
if $X$ admits two such models $\cX$ and $\cX'$,
the class of the special fibres $\cX_0$ and $\cX'_0$
in
$\kvar_{\rm loc}$ are equal, which may be seen as an analogue
of Batyrev's result on birational Calabi-Yau varieties
\cite{BCY}.

The paper is organized as follows. Section \ref{sec2} is devoted to preliminaries
on formal schemes, the Greenberg functor and weak N{\'e}ron models.
In the following section, we review the results on motivic integration
on formal schemes obtained by the second named author in \cite{sebag}
which are needed in the present work. 
We are then able in section \ref{sec4} to construct a  motivic integration 
on smooth rigid varieties and to prove the main results which were mentionned
in the present introduction.
Finally, in section \ref{sec5}, guided by the analogy with arc spaces,
we formulate an analogue of the Nash problem, which is about  the
relation between essential (\textit{i.e.}
appearing in every resolution) components
of resolutions of a singular variety and irreducible components
of spaces of truncated arcs on the variety,
for formal $R$-schemes with smooth generic fibre. 
In this context, analogy suggests there might be some
relation between essential components of weak N{\'e}ron models 
of a given
formal $R$-scheme $\cX$ with smooth generic fibre
and irreducible 
components of the truncation $\pi_n (\Gr (\cX))$ of its Greenberg
space for $n \gg 0$. As a very first step in that direction
we compute 
the dimension of the contribution of a given irreducible 
component to the truncation.

\section{Preliminaries on formal schemes and
Greenberg functor}\label{sec2}
\subsection{Formal schemes}In this paper $R$ will denote a complete discrete valuation 
ring with residue field $k$ and fraction field $K$. We shall assume $k$ is
perfect. We shall fix once for all an uniformizing parameter $\varpi$
and we shall set $R_n := R / (\varpi)^{n +1}$, for $n \geq 0$.
In the whole paper,
by  a formal $R$-scheme, we shall always mean a quasi-compact, separated,
locally topologically
of finite type formal $R$-scheme, in the sense of
\S\kern .15em 10 of \cite{EGA}.
A formal $R$-scheme is a locally ringed space $(\cX, \cO_{\cX})$
in topological $R$-algebras.
It is equivalent to the data,
for every $n \geq 0$, of the
$R_n$-scheme $X_n = (\cX, \cO_{\cX} \otimes_R R_n)$.
The $k$-scheme $X_0$ is called
the special fibre of $\cX$. As a topological space 
$\cX$ is isomorphic to $X_0$ and $\cO_{\cX} = \limproj \cO_{X_n}$.
We have $X_n = X_{n + 1} \otimes_{R_{n +1}} R_n$
and $\cX$
is canonically isomorphic to the inductive limit of the schemes
$X_n$ in the category of formal schemes.
Locally $\cX$ is isomorphic to an affine 
formal $R$-scheme of the form ${\Specf}  A$ with $A$ an
$R$-algebra topologically of finite type, \textit{i.e.}  a quotient of a
restricted formal series
algebra 
$R \{T_1, \dots, T_m\}$. If $\cY$ and $\cX$
are $R$-formal schemes,
we denote by ${\rm Hom}_R (\cY, \cX)$ the set of morphisms
of formal $R$-schemes $\cY \rightarrow \cX$, \textit{i.e.}
morphisms between the underlying
locally topologically ringed spaces over $R$ 
(cf.
\S\kern .15em 10 of \cite{EGA}). 
It follows from
Proposition 10.6.9 of \cite{EGA}, that the canonical morphism
${\rm Hom}_R (\cY, \cX) \rightarrow \limproj {\rm Hom}_{R_n} (\cY_n, \cX_n)$
is a bijection.

If $k$ is a field, by a variety over $k$ we
mean
a separated reduced scheme of finite type over $k$.

\subsection{Extensions}Let $A$ be a $k$-algebra.
We set
$L (A) = A$  when $R$ is a ring of equal characteristic  
and $L (A) = W (A)$, the ring of Witt vectors,
when
$R$ is a ring of unequal characteristic and we denote
by $R_A$ the ring $R_A := R \otimes_{L (k)} L(A)$.
When $F$ is a field
containing
$k$, we denote by $K_F$
the field of fractions of $R_F$. When the  field $F$
is perfect, the ring $R_F$ is a discrete valuation ring and, furthermore,
the uniformizing parameter
$\varpi$ in $R$ induces an uniformizing parameter in $R_F$.
Hence, since $k$ is assumed to be perfect,
the extension $R \rightarrow R_F$
has ramification index 1 in the terminology of
\S\kern .15em 3.6 of \cite{neron}.

\subsection{The Greenberg Functor}We shall recall some material from
\cite{g1}
and \S\kern .15em 9.6 of \cite{neron}.
Let us remark, that, when $R$ is a ring of equal characteristic,
we can view $R_n$ as the set of $k$-valued points of some affine space
$\AA^m_k$ which we shall denote by $\cR_n$, in a way compatible
with the $k$-algebra structure.
When $R$ is a ring of unequal characteristic,
$R_n$ can no longer be viewed as a $k$-algebra. However, using Witt vectors,
we may still
interpret $R_n$ as the set of $k$-valued points of a ring scheme
$\cR_n$, which, as a $k$-scheme, is isomorphic 
to some affine space $\AA^m_k$. Remark we have canonical morphisms
$\cR_{n + 1} \rightarrow \cR_n$.

Now, for every $n \geq 0$,
we consider the functor $h_n^{\ast}$ which to a $k$-scheme $T$ associates
the locally ringed space $h_n^{\ast} (T)$
which has $T$ as underlying topological space and ${\cal Hom}_k (T, \cR_n)$
as structure sheaf.
In particular, for any $k$-algebra $A$,
$$
h_n^{\ast} (A) = \Spec (R_n \otimes_{L (k)} L (A)).
$$
Taking $A = k$, we see that $h_n^{\ast} T$ is a locally ringed 
space over $\Spec R_n$.

By a fundamental result of Greenberg \cite{g1} 
(which in
the equal characteristic case amounts to Weil restriction of scalars),
for $R_n$-schemes $X_n$, locally of finite type,
the functor 
$$
T \longmapsto {\rm Hom}_{R_n} (h_n^{\ast} (T), X_n)
$$
from the category of $k$-schemes to the category of sets
is represented by a $k$-scheme $\Gr_n (X_n)$
which is locally of finite type.
Hence, for any $k$-algebra $A$,
$$
\Gr_n (X_n) (A) = X_n (R_n \otimes_{L (k)} L (A)),
$$
and, in particular, setting $A =k$, we have
$$
\Gr_n (X_n) (k) = X_n (R_n).
$$

Among basic properties the Greenberg functor respects
closed immersions, open imersions, fibred products, smooth and {\'e}tale morphisms
and also sends affines to affines.

Now let us consider again
a formal $R$-scheme $\cX$.
The canonical adjunction morphism
$
h_{n +1}^{\ast} (\Gr_{n +1} (X_{n +1})) \rightarrow X_{n +1}
$
gives rise, by tensoring with $R_n$, to a canonical morphism
of of $R_n$-schemes
$
h_{n}^{\ast} (\Gr_{n +1} (X_{n +1})) \rightarrow X_{n},
$
from which one derives, again by adjunction, 
a canonical morphism of
$k$-schemes
$$
\theta_n : \Gr_{n +1} (X_{n +1}) \longrightarrow \Gr_{n} (X_{n}).
$$
In this way we attach to the
formal scheme $\cX$ a projective system
$(\Gr_{n} (X_{n}))_{n \in \NN}$ of $k$-schemes. The 
transition morphisms
$\theta_n$ being affine, the projective
limit
$$
\Gr (\cX) := \limproj \Gr_{n} (X_{n})
$$
exists
in the category of $k$-schemes.

Let $T$ be a 
$k$-scheme. We denote by $h^{\ast} (T)$
the locally ringed space
which has $T$ as underlying topological space and 
$\limproj {\cal Hom}_k (T, \cR_n)$
as structure sheaf.
It is a locally ringed space over $\Specf R$ which identifies
with the projective limit of the spaces
$h_n^{\ast} (T)$ in the category of 
locally ringed spaces.
Furthermore one checks, similarly as in
Proposition 10.6.9 of \cite{EGA}, that the canonical morphism
${\rm Hom}_R (h^{\ast} (T), \cX) \rightarrow \limproj {\rm Hom}_{R_n}
(h_n^{\ast} (T), \cX_n)$
is a bijection for every formal $R$-scheme $\cX$.

Putting everything together we get the following:
\begin{prop}Let $\cX$ be a 
quasi-compact, 
locally topologically
of finite type formal $R$-scheme. The functor 
$$
T \longmapsto {\rm Hom}_{R} (h^{\ast} (T), \cX)
$$
from the category of $k$-schemes to the category of sets
is represented by the  $k$-scheme $\Gr (\cX)$. \qed
\end{prop}

In particular, for every field $F$ containing
$k$,
there are  canonical bijections
$$
\Gr (\cX) (F) \simeq {\rm Hom}_R (\Specf R_F, \cX) \simeq
\cX (R_F).
$$

One should note that, in general,
$\Gr (\cX)$ is not of finite type,
even if  $\cX$  is a quasi-compact, 
topologically
of finite type formal $R$-scheme.

In this paper, we shall always consider the schemes
$\Gr_{n} (X_{n})$ and $\Gr (\cX)$
with their reduced structure.

Sometimes, by abuse of notation,  we shall write
$\Gr_n (\cX)$ for $\Gr_n (X_n)$.

\begin{prop}\label{glue}
\begin{enumerate}
\item[(1)]The functor $\Gr$ respects
open and closed immersions, fibre products, and sends
affine formal $R$-schemes
to affine $k$-schemes.
\item[(2)]
Let $\cX$ be a formal quasi-compact and separated
$R$-scheme and let $(\cO_i)_{i \in J}$
be a finite covering by formal open subschemes.
There are canonical isomorphisms $\Gr (\cO_i \cap \cO_j) \simeq
\Gr (\cO_i) \cap \Gr (\cO_j)$ and 
the scheme $\Gr (\cX)$ is canonically isomorphic
to the scheme obtained by glueing the schemes
$\Gr (\cO_i)$.
\end{enumerate}
\end{prop}

\begin{proof}Assertion (1) for the functor $\Gr_n$
is proved in \cite{g1}
and \cite{neron}, and follows for
$\Gr$ by taking projective limits.
Assertion (2) follows from (1) and the universal property
defining $\Gr$.
\end{proof}

\begin{remark}\label{arcspace}Assume we are in
the equal characteristic case, \textit{i.e.} $R = k[[\varpi]]$.
For $X$ an algebraic  variety over $k$, we can consider 
the formal $R$-scheme $X \widehat \otimes R$
obtained by base change
and completion.
We have
canonical isomorphisms
$\Gr (X \widehat \otimes R) \simeq \cL (X)$
and
$\Gr_n (X \otimes R_n) \simeq \cL_n (X)$,
where 
$\cL (X)$ and $\cL_n (X)$
are the arc spaces considered in \cite{arcs}.
\end{remark}

\subsection{Smoothness}\label{smooth}
Let us recall the definition of smoothness
for morphisms of formal $R$-schemes. A morphism
$f : \cX \rightarrow \cY$ of formal $R$-schemes is smooth at a point $x$ of
$X_0$ of relative dimension $d$ if it is flat at $x$
and the induced morphism $f_0 : X_0 \rightarrow Y_0$ is smooth
at $x$ of relative dimension $d$.
An equivalent condition (cf. Lemma 1.2 of \cite{rigid2}) is that 
for every $n$ in $\NN$ the induced morphism
$f_n : X_n \rightarrow Y_n$ is smooth
at $x$ of relative dimension $d$.
The morphism $f$ is smooth if it is smooth at every point of $X_0$.
The  formal $R$-scheme $\cX$ is smooth at a point $x$ of
$X_0$ if the structural morphism is smooth at $x$.

Let $\cX$ be a flat formal $R$-scheme of relative dimension $d$.
We denote by $\cX_{\rm sing}$ the closed 
formal subscheme of $\cX$ defined
by the radical of the Fitting ideal sheaf ${\rm Fitt}_d \Omega_{\cX / R}$.
The  formal $R$-scheme $\cX$ is smooth at a point $x$ of
$X_0$ (resp. is smooth)
if and only if $x$ is not in $\cX_{\rm sing}$ (resp. $\cX_{\rm sing}$
is empty).

\subsection{Greenberg's Theorem}The following statement, which is an adaptation of a result of
Schappacher \cite{scha}, is an 
analogue of Greenberg's Theorem
\cite{g2} in the framework of formal schemes.
We refer to \cite{sebag} for a more detailled exposition.

\begin{theorem}\label{gree}Let $R$ be a complete discrete valuation 
ring and let $\cX$ be a formal $R$-scheme. For every $n \geq 0$, there exists
an integer $\gamma_{\cX} (n) \geq n$ such that, for every 
field $F$ containing $k$, and every $x$ in
$\cX (R_F/ \varpi^{\gamma_{\cX} (n) +1})$, the image of $x$ in 
$\cX (R_F/\varpi^{n +1})$
may be lifted to a point in
$\cX (R_F)$.
\end{theorem}

\begin{remark}\label{gf}The function $n \mapsto \gamma_{\cX} (n)$
is called the Greenberg function of $\cX$.
\end{remark}

\subsection{Rigid spaces}For  $\cal X$ a flat formal $R$-scheme
we shall denote by $\cX_{K}$ its generic fibre in the sense of Raynaud
\cite{table}. By Raynaud's Theorem \cite{table}, \cite{rigid1},
the functor 
$\cX \mapsto \cX_{K}$ induces an equivalence between the 
localization of the category
of quasi-compact flat formal $R$-schemes by admissible
formal blowing-ups, and the category of rigid $K$-spaces
which are quasi-compact and quasi-separated.
Furthemore, $\cX$ is separated if and only if
$\cX_{K}$ is separated (cf. Proposition 4.7 of \cite{rigid1}).
Recall that for the blowing-up of an ideal sheaf $\cI$ to be admissible
means that $\cI$ contains some power of the
uniformizing parameter $\varpi$.
In the paper all rigid $K$-spaces will be assumed to be
quasi-compact and separated.

\subsection{Weak N{\'e}ron models}We shall denote by $R^{\rm sh}$
a strict henselization of 
$R$ and by 
$K^{\rm sh}$ its field of fractions.
\begin{definition}Let $X$ be a smooth rigid $K$-variety.
A weak formal N{\'e}ron model\footnote{We follow here the
terminology of \cite{bs} which is 
somewhat different from that of \cite{neron}.} of $X$ is a smooth formal $R$-scheme $\cU$,
whose generic fibre $\cU_K$ is an open rigid subspace of $X_K$, and which has
the property that the canonical map
$\cU (R^{\rm sh}) \rightarrow X (K^{\rm sh})$
is bijective.
\end{definition}

The construction of weak N{\'e}ron models using N{\'e}ron's 
smoothening process presented
in \cite{neron}, carries over almost literally from $R$-schemes to formal
$R$-schemes, and gives, as explained in \cite{bs},
the following result:

\begin{theorem}\label{wn}Let $\cX$ be a quasi-compact formal $R$-scheme,
whose generic fibre $\cX_K$ is smooth.
Then there exists a morphism of formal $R$-schemes $\cX'\rightarrow \cX$,
which is the composition of a sequence of formal blowing-ups with centers in
the corresponding
special fibres, such that every
$R^{\rm sh}$-valued point of $\cX$ factors through the smooth locus of $\cX'$.
\end{theorem}

One deduces the following omnibus statement:
\begin{prop}Let $X$ be a smooth quasi-compact and separated rigid
$K$-space and let $\cX$ be a formal $R$-model  of $X$, \textit{i.e} a
quasi-compact formal $R$-scheme $\cX$ with generic fibre $X$. Then
there exists
a weak formal N{\'e}ron model $\cU$ of $X$ which dominates
$\cX$ and which is quasi-compact.
Furthermore the canonical map $\cU (R^{\rm sh}) \rightarrow 
\cX (R^{\rm sh})$ is a bijection and 
for every perfect field $F$ containing $k$,
the formal $R_F$-scheme
$\cU \otimes_R R_F$ is a weak N{\'e}ron model of 
the rigid $K_F$-space $X \otimes_K K_F$.
In particular, the morphism $\cU \rightarrow \cX$ 
induces a bijection between points of
$\Gr (\cU)$ and $\Gr (\cX)$.
\end{prop}

\begin{proof}We choose a formal model $\cX$  of $X$ such that we are
in the situation of Theorem \ref{wn}. The smooth locus $\cU$
of $\cX'$ is quasi-compact and is a weak N{\'e}ron model of
$\cX_K$, since, by \cite{bs} 2.2 (ii), every $K^{\rm sh}$-valued point
of
$\cX_K$ extends uniquely to a
$R^{\rm sh}$-valued point of $\cX$.
Also, it follows from Corollary 6 of \S\kern .15em 3.6
\cite{neron} that,
if $\cU$ is a weak N{\'e}ron model of
the rigid $K$-space $X$, then for every field $F$ containing $k$,
the formal $R_F$-scheme
$\cU \otimes_R R_F$ is a weak N{\'e}ron model of 
the rigid $K_F$-space $X \otimes_K K_F$.
\end{proof}

\section{Motivic integration on formal schemes}\label{sec3}
The material in this section is borrowed from \cite{sebag}
to which we shall refer for details.

\subsection{Truncation}For $\cX$ a formal $R$-scheme, we shall denote
by $\pi_{n, \cX}$ or $\pi_n$ the canonical projection
$\Gr (\cX) \rightarrow \Gr_n (X_n)$, for $n$ in $\NN$.

Let us first state the following corollary
of Theorem \ref{gree}:

\begin{prop}\label{chev}Let $\cX$ be a formal $R$-scheme.
The image $\pi_n (\Gr (\cX))$ of 
$\Gr (\cX)$ in $\Gr_n (X_n)$ is a constructible subset of
$\Gr_n (X_n)$. More generally, if $C$ is a constructible
subset of $\Gr_m (X_m)$, 
$\pi_n (\pi_m^{-1} (C))$  is a constructible subset of
$\Gr_n (X_n)$, for every $n \geq 0$.
\end{prop}

\begin{proof}Indeed, it follows from
Theorem \ref{gree} that
$\pi_n (\Gr (\cX))$ is equal to the image of
$\Gr_{\gamma (n)} (X_{\gamma (n)})$
in
$\Gr_n (X_n)$. The morphism  
$\Gr_{\gamma (n)} (X_{\gamma (n)}) \rightarrow \Gr_n (X_n)$
being of finite type, first the statement follows from
Chevalley's Theorem. For the second statement, one may assume $m = n$, and the proof proceeds as before.
\end{proof}

\begin{prop}\label{smoothtr}
Let $\cal X$ be a smooth
formal separated $R$-scheme (quasi-compact, 
locally topologically
of finite type over $R$), of relative dimension $d$.
\begin{enumerate}
\item[(1)]For every $n$, the morphism
$\pi_n : \Gr (\cX) \rightarrow \Gr_n (X_n)$ is surjective.
\item[(2)]For every $n$ and $m$ in $\NN$, the canonical
projection
$\Gr_{n +m} (X_{n +m}) \rightarrow \Gr_n (X_n)$
is a locally trivial fibration for the Zariski topology with fibre
$\AA_k^{dm}$.
\end{enumerate}
\end{prop}

We say that
a map
$\pi : A \rightarrow B$ is a  piecewise morphism if
there exists a finite partition of
the domain of $\pi$ into locally closed subvarieties of $X$
such that the restriction of $\pi$ to any of these
subvarieties is a morphism of
schemes.

\subsection{Away from the singular locus}Let $\cX$ be a formal $R$-scheme and
consider its singular locus 
$\cX_{\rm sing}$ defined  in \ref{smooth}. For any integer $e \geq 0$, we view 
$\Gr_{e} (\cX_{{\rm sing}, e})$ as contained in 
$\Gr_{e} (\cX)$ and we set
$$
\Gr^{(e)} (\cX) := \Gr (\cX) \setminus \pi_e^{-1} (\Gr_{e} (\cX_{{\rm sing}, e})).
$$

We say that
a map 
$\pi : A \rightarrow B$ between $k$-constructible sets $A$ and $B$
is a
piecewise trivial fibration with fibre
$F$, if there exists a finite partition of $B$ in subsets $S$ which are
locally closed
in $Y$ such that $\pi^{- 1} (S)$ is locally closed in $X$ and
isomorphic, as
a variety over $k$, to $S \times F$, with $\pi$
corresponding under the isomorphism to the projection
$S \times F \rightarrow S$. We say that the map $\pi$ is
a piecewise trivial fibration over some constructible subset $C$ of
$B$,
if the restriction of $\pi$ to $\pi^{- 1} (C)$ is a piecewise 
trivial fibration onto $C$.

\begin{prop}\label{ptf}Let
$\cX$ a flat formal $R$-scheme of relative dimension $d$.
There exists an integer $c
\geq 1$
such that, for every integers $e$ and $n$ in $\NN$ such that
$n \geq c e$,
the projection
$$
\pi_{n+1} (\Gr (\cX)) \longrightarrow \pi_{n} (\Gr (\cX))
$$
is a piecewise trivial fibration over $\pi_{n} (\Gr^{(e)} (\cX))$
with fibre $\AA_k^d$.
\end{prop}

\subsection{Dimension estimates}

\begin{lem}\label{dimgr}
Let $\cX$ be a formal $R$-scheme whose generic fibre $\cX_K$ is of dimension
$\leq d$.
Then
\begin{enumerate}
\item[(1)]For every $n$ in $\NN$, $${\rm dim} \, \pi_n (\Gr (\cX))
\leq (n + 1) d.$$
\item[(2)]For $m \geq n$, the fibres
of the projection
$\pi_m (\Gr (\cX)) \rightarrow \pi_n (\Gr (\cX))$
are of dimension $\leq (m - n) d$.
\end{enumerate}
\end{lem}

\begin{lem}\label{dim2}
Let $\cX$ be a formal $R$-scheme whose generic fibre $\cX_K$ is of dimension
$d$. Let $\cS$ be a closed 
formal $R$-subscheme of $\cX$ such that $\cS_K$ is of dimension
$< d$.
Then, for all integers $n$, $i$ and $\ell$ such that
$n \geq i \geq \gamma_{\cS} (\ell)$,
where $\gamma_{\cS}$ is the Greenberg function of $\cS$ defined in \ref{gf},
$\pi_{n, \cX} (\pi_{i, \cX}^{-1} \Gr_i (\cS))$
is of dimension
$\leq (n + 1) d - \ell -1$.
\end{lem}

\subsection{Grothendieck rings}\label{GK}
Let $k$ be a field.
We denote by $\kvar$ the abelian group
generated by symbols $[S]$, for $S$ a variety over
$k$, with the relations $[S] = [S']$ if $S$ and $S'$ are
isomorphic and $[S] = [S'] + [S \setminus S']$
if $S'$ is closed in $S$. There is a natural ring structure
on $\kvar$, the product being
induced by the cartesian product of varieties, and
to any constructible set $S$ in some variety  one naturally associates
a class $[S]$ in $\kvar$.
We denote
by $\kvar_{\rm loc}$ the localisation 
$\kvar_{\rm loc} :=  \kvar [\LL^{-1}]$ with
$\LL := [\AA^1_{k}]$. We denote by 
$F^m \kvar_{\rm loc}$ the subgroup generated by
$[S] \LL^{- i}$ with ${\rm dim} \, S - i \leq -m$, and by
$\widehat{\kvar}$ the completion
of $\kvar_{\rm loc}$ with respect to
the filtration $F^{\cdot}$\footnote{It is still unknown whether
the filtration $F^{\cdot}$ is separated or not.}. We
will also denote by $F^{\cdot}$ the filtration induced on
$\widehat{\kvar}$. 
We denote by $\overline \kvar_{\rm loc}$ the image of 
$\kvar_{\rm loc}$ in $\widehat \kvar$.
We put on the ring 
$\widehat{\kvar}$
a structure of non-archimedean ring by setting
$||a|| := 2^{-n}$,
where $n$ is the largest $n$ such that $a \in F^{n} \widehat \kvar$, for $a \not=0$ and $||0|| = 0$.

\subsection{Cylinders}Let $\cX$ be a formal $R$-scheme.
A subset $A$ of $\Gr (\cX)$ is cylindrical of level $n \geq 0$ if
$A = \pi_n^{- 1} (C)$
with $C$ a constructible subset of 
$\Gr_n (\cX)$.
We denote by $\CC_{\cX}$ the set of cylindrical subsets
of $\Gr (\cX)$ of some level.
Let us remark that
$\CC_{\cX}$ is a boolean algebra, \textit{i.e.} contains
$\Gr (\cX)$, $\emptyset$, and is  stable by finite intersection, finite union,
and by taking complements.
It follows from Proposition \ref{chev}, that if $A$ is cylindrical of 
some level, then $\pi_n (A)$ is constructible for every $n \geq 0$.

A basic finiteness property of cylinders is the following:
\begin{lem}\label{finiteness}
Let $A_i$, $i \in I$, be an enumerable family of
cylindrical subsets of $\Gr (\cX)$. If $A := \cup_{i \in I} A_i$ is
also a cylinder, then there exists a finite subset $J$ of $I$ such that
$A := \cup_{i \in J} A_i$.
\end{lem}

\begin{proof}Since $\Gr (\cX)$ is quasi-compact, this follows from 
Th{\'e}or{\`e}me 7.2.5 of \cite{EGA}.
\end{proof}

\subsection{Motivic measure for cylinders}Let
$\cX$ a flat formal $R$-scheme of relative dimension $d$.
Let $A$ be a cylinder of $\Gr (\cX)$. We shall say
$A$ is stable of level $n$ if
it is cylindrical of level $n$ and if, for every $m \geq n$,
the morphism 
$$
\pi_{m+1} (\Gr (\cX)) \longrightarrow \pi_{m} (\Gr (\cX))
$$
is a piecewise trivial fibration over $\pi_{n} (A)$
with fibre $\AA_k^d$.
We denote by $\CC_{0, \cX}$ the set of stable
cylindrical subsets
of $\Gr (\cX)$ of some level.

It follows from Proposition \ref{smoothtr} that every cylinder in
$\Gr (\cX)$ is stable
when $\cX$ is smooth. When $\cX$ is no longer assumed to be smooth,
$\CC_{0, \cX}$ is in general not a boolean algebra, but is an ideal of
$\CC_{\cX}$: $\CC_{0, \cX}$ contains 
$\emptyset$, is stable by finite union, and the intersection
of an element in $\CC_{\cX}$ with an element of 
$\CC_{0, \cX}$ belongs to $\CC_{0, \cX}$. In general $\Gr (\cX)$ is not
stable, but it
follows from Proposition \ref{ptf} that
$\Gr^{(e)} (\cX)$ is a stable cylinder of $\Gr (\cX)$, for every $e \geq 0$.

From  first principles, one proves (cf. \cite{arcs}, \cite{sebag}):

\begin{def-prop}There is a unique additive morphism
$$\tilde \mu : \CC_{0, \cX} \longrightarrow \kvar_{\rm loc}$$
such that $\tilde \mu (A) = [\pi_n (A)] \, \LL^{- (n +1) d}$,
when $A$ is a stable cylinder of level $n$.
Furthermore the measure $\tilde \mu$ is $\sigma$-additive.
\end{def-prop}

One deduces from Proposition \ref{dimgr} and Proposition \ref{dim2}, cf.
\cite{arcs}, \cite{sebag}, the following:

\begin{prop}\label{conv}
\begin{enumerate}
\item[(1)]For any  cylinder $A$ in 
$\CC_{\cX}$, the limit
$$
\mu (A) := \lim_{e \rightarrow \infty} \tilde \mu (A \cap \Gr^{(e)} (\cX))
$$
exists in $\widehat \kvar$.
\item[(2)]If $A$ belongs to $\CC_{0, \cX}$, then $\mu (A)$ coincide with the
image of $\tilde \mu (A)$ in $\widehat \kvar$.
\item[(3)]The measure $A \mapsto \mu (A)$ is $\sigma$-additive.
\item[(4)]For $A$ and $B$ in $\CC_{\cX}$,
$|| \mu (A \cup B)|| \leq \max (||\mu (A)||, ||\mu (B)||)$. If $A \subset B$,
$|| \mu (A)|| \leq || \mu (B)||$.
\end{enumerate}
\end{prop}

\subsection{Measurable subsets of $\Gr (\cX)$}
For $A$ and $B$ subsets of the same set, we use the notation
$A \triangle B$ for $A \cup B \setminus A \cap B$.
\begin{definition}We say that a subset $A$ of $\Gr (\cX)$ is
{\emph{measurable}}
if, for every positive real number $\varepsilon$, there exists an
$\varepsilon$-cylindrical approximation, \textit{i.e.} a
sequence of
cylindrical subsets
$A_{i} (\varepsilon)$, $i \in \NN$,
such that
$$
\Bigl(A \triangle
A_{0} (\varepsilon)
\Bigr)
\subset \bigcup_{i \geq 1} A_{i} (\varepsilon),
$$
and 
$||\mu (A_{i} (\varepsilon))|| \leq \varepsilon$
for all $i \geq 1$. We say that $A$ is {\emph{strongly measurable}}
if moreover we can take $A_{0} (\varepsilon) \subset A$.
\end{definition}

\begin{theorem}If $A$ is a measurable subset of $\Gr (\cX)$,
then 
$$\mu (A) :=
\lim_{\varepsilon \rightarrow 0}  \mu (A_{0} (\varepsilon))$$
exists in $\widehat \kvar$ and is independent of the choice of the
sequences
$A_{i} (\varepsilon)$, $i \in \NN$.
\end{theorem}

For $A$ a measurable subset of $\Gr (\cX)$, we shall call
$\mu (A)$ the motivic measure of $A$.
We denote by $\DD_{\cX}$ the set of 
measurable subsets of $\Gr (\cX)$.

One should remark that obviously
$\CC_{\cX}$ is contained in $\DD_{\cX}$.

\begin{prop}\begin{enumerate}
\item[(1)]$\DD_{\cX}$ is a boolean algebra.
\item[(2)]If $A_{i}$, $i \in \NN$, is a sequence
of measurable subsets
of $\Gr (\cX)$ with
$\lim_{i \rightarrow \infty} ||\mu (A_{i})|| = 0$, then
$\cup_{i \in \NN} A_{i}$ is measurable.
\item[(3)] Let $A_{i}$, $i \in \NN$, be a family of measurable
subsets of $\Gr (\cX)$. Assume the sets $A_{i}$ are mutually disjoint
and that $A := \cup_{i \in \NN}A_{i}$ is measurable.
Then $\sum_{i \in \NN} \mu (A_{i})$
converges in $\widehat \kvar$ to $\mu (A)$.
\item[(4)] If $A$ and $B$ are measurable
subsets of $\Gr (\cX)$ and if $A \subset B$, then
$||\mu (A)|| \leq ||\mu (B)||$.
\end{enumerate}
\end{prop}

\begin{remark}\label{comp}In the situation of Remark \ref{arcspace},
one can check that the notion of cylinders, stable cylinders, measurable
subsets of $\Gr (X \widehat \otimes R)$ coincides with the analogous
notions
introduced in \cite{MK} for subsets of $\cL (X)$.
\end{remark}

\subsection{Order of the Jacobian ideal}
Let $h : \cY \rightarrow \cX$ be a morphism of
flat formal $R$-schemes of relative dimension $d$.

Let $y$ be a point of $\Gr (\cY) \setminus \Gr (\cY_{\rm sing})$
defined over some field extension $F$ of $k$.
We denote by $\varphi : \Specf R_F \rightarrow \cY$
the corresponding morphism of formal $R$-schemes.
We define 
$\ord_{\varpi} (\Jac_h) (y)$, 
the order of the Jacobian ideal of $h$ at $y$, as follows.

From the natural morphism
$h^{\ast} \Omega_{\cX | R} \rightarrow \Omega_{\cY | R}$,
one deduces, by taking the $d$-th exterior power,
a morphism
$h^{\ast} \Omega^d_{\cX | R} \rightarrow \Omega^d_{\cY | R}$,
hence a morphism
$$(\varphi^{\ast}h^{\ast} \Omega^d_{\cX | R}) /(\text{torsion})
\longrightarrow (\varphi^{\ast}\Omega^d_{\cY | R})/(\text{torsion}).$$
Since  $L:= (\varphi^{\ast}\Omega^d_{\cY | R})/(\text{torsion})$
is a free $\cO_{R_F}$-module of rank 1, it follows from the structure
theorem for finite type modules over principal domains,
that the image of
$M := (\varphi^{\ast}h^{\ast} \Omega^d_{\cX | R}) /(\text{torsion})$
in $L$ is either $0$, in which case we set
$\ord_{\varpi} (\Jac_h) (y) =\infty$,
or $\varpi^n L$, for some $n \in \NN$, in which case we set
$\ord_{\varpi} (\Jac_h) (y) = n$.

\subsection{The change of variable formula}
If $h :\cY \rightarrow \cX$ is a morphism of formal
$R$-schemes, we still write $h$ for the corresponding
morphism $\Gr (\cY) \rightarrow \Gr (\cX)$.

\medskip

The following lemmas are basic geometric ingredients
in the proof of the change of variable formula \ref{cv} and \ref{cvtilde}.

\begin{lem}\label{cv1}Let $h : \cY \rightarrow \cX$ be a
morphism between
flat formal $R$-schemes of relative dimension $d$.
We assume $\cY$ is smooth.
For $e$ and 
$e'$ in $\NN$ we
set
$$
\Delta_{e, e'} := 
\Bigl\{ \varphi \in \Gr (\cY) \Bigm \vert 
\ord_{\varpi} (\Jac_h) (y) = e \quad \text{\rm and} \quad h(\varphi)
\in \Gr^{(e')} (\cX) \Bigr\}.
$$
Then, there exists $c$ in $\NN$ such that,
for every $n \geq 2e$, $n \geq e + c e'$,
for every $\varphi$ in $\Delta_{e, e'}$ and every $x$ in $\Gr (\cX)$ such that
$\pi_n (h (\varphi)) = \pi_n (x)$, there exists $y$ in
$\Gr (\cY)$ such that $h (y) = x$ and
$\pi_{n - e} (\varphi) = \pi_{n - e} (y)$.
\end{lem}

\begin{lem}\label{cv2}Let $h : \cY \rightarrow \cX$ be a
morphism between
flat formal $R$-schemes of relative dimension $d$.
We assume $\cY$ is smooth.
Let $B$ be a cylinder in $\Gr (\cY)$ and set $A = h (B)$.
Assume 
$\ord_{\varpi} (\Jac_h)$ is constant with value $e < \infty$ on $B$
and that $A \subset \Gr^{(e')} (\cX)$ for some $e' \geq 0$.
Then $A$ is a cylinder. Furthermore, if the restriction of $h$ to $B$
is injective, then for $n \gg 0$, the following holds:
\begin{enumerate}
\item[(1)]If $\varphi$ and $\varphi'$ belong to $B$ and
$\pi_n (h (\varphi)) = \pi_n (h (\varphi'))$,
then $\pi_{n - e} (\varphi) = \pi_{n - e} (\varphi')$.
\item[(2)]The morphism $\pi_n (B) \rightarrow \pi_n (A) $ induced by $h$
is a
piecewise trivial fibration with fibre $\AA_k^e$.
\end{enumerate}
\end{lem}

For a measurable subset $A $ of $\Gr (\cX)$ and a function $\alpha : A
\rightarrow
\ZZ \cup \{\infty\}$, we say that $\LL^{-\alpha}$ is {{integrable}}
or that $\alpha$ is {{exponentially integrable}}
if the fibres of $\alpha$ are measurable and if the motivic integral
$$
\int_{A} \LL^{- \alpha} d\mu
:= 
\sum_{n \in \ZZ} \mu (A \cap \alpha^{-1} (n)) \LL^{- n}
$$
converges in $\widehat \kvar$.
 
When all the fibres $A \cap \alpha^{-1} (n)$ are stable cylinders 
and $\alpha$ takes only a finite number of values on $A$,
it is not necessary to go to the completion of
$\kvar_{\rm loc}$ and one may define directly
$$
\int_{A} \LL^{- \alpha} d\tilde \mu
:= 
\sum_{n \in \ZZ} \tilde \mu (A \cap \alpha^{-1} (n)) \LL^{- n}
$$
in $\kvar_{\rm loc}$.

\begin{theorem}\label{cv}Let $h : \cY \rightarrow \cX$ be a
morphism between
flat formal $R$-schemes of relative dimension $d$.
We assume $\cY$ is smooth.
Let $B$ be a strongly measurable subset of $\Gr (\cY)$.
Assume $h$ induces a bijection between $B$ and $A := h(B)$.
Then, for every exponentially integrable function
$\alpha : A \rightarrow \ZZ \cup \infty$,
the function
$\alpha \circ h + \ord_{\varpi}(\Jac_h)$ is 
exponentially integrable on $B$ and
$$
\int_A \LL^{- \alpha} d\mu =
\int_B \LL^{-\alpha \circ h - \ord_{\varpi}(\Jac_h)} d\mu.
$$
\end{theorem}

We shall also need the following variant of Theorem \ref{cv}.
\begin{theorem}\label{cvtilde}Let $h : \cY \rightarrow \cX$ be a
morphism between
flat formal $R$-schemes of relative dimension $d$.
We assume $\cY$ and $\cX_K$ are  smooth and that 
the morphism $h_K : \cY_K \rightarrow \cX_K$ induced by
$h$ 
is {\'e}tale (cf. \cite{rigid3}).
Let $B$ be a cylinder in $\Gr (\cY)$.
Assume $h$ induces a bijection between $B$ and $A := h(B)$ and that
$A$ is a stable cylinder of $\Gr (\cX)$.
Then
the fibres $B \cap \ord_{\varpi}(\Jac_h)^{-1} (n)$ are stable cylinders,
$\ord_{\varpi}(\Jac_h)^{-1} (n)$ takes only a finite number of values on $B$,
and 
$$
\int_A d \tilde \mu =
\int_B \LL^{- \ord_{\varpi}(\Jac_h)} d \tilde\mu
$$
in $\kvar_{\rm loc}$.
\end{theorem}

\section{Integration
on smooth rigid varieties}\label{sec4}

\subsection{Order of differential forms}
Let
$\cX$ be a flat formal $R$-scheme equidimensional
of relative dimension $d$.
Consider a differential form $\omega$ in $\Omega^d_{\cX | R} (\cX)$.
Let $x$ be a point of $\Gr (\cX) \setminus \Gr (\cX_{\rm sing})$
defined over some field extension $F$ of $k$.
We denote by $\varphi : \Specf R_F \rightarrow \cX$
the corresponding morphism of formal $R$-schemes.
Since  $L:= (\varphi^{\ast}\Omega^d_{\cX | R})/(\text{torsion})$
is a free $\cO_{R_F}$-module of rank 1, it follows from the structure
theorem for finite type modules over principal domains,
that its submodule
$M$ generated by $\varphi^{\ast} \omega$
is either $0$, in which case we set
$\ord_{\varpi} (\omega) (x) = \infty$,
or $\varpi^n L$, for some $n \in \NN$, in which case we set
$\ord_{\varpi} (\omega) (x) = n$.

Since there is a canonical isomorphism
$\Omega^d_{\cX_{K}} (\cX_K) \simeq \Omega^d_{\cX | R} (\cX) \otimes_R K$
(cf. Proposition 1.5 of \cite{rigid3}),
if $\omega$ is in $\Omega^d_{\cX_{K}} (\cX_K)$, 
we write $\omega = \varpi^{-n}
\tilde \omega$, with
$\tilde \omega$ in $\Omega^d_{\cX | R} (\cX)$ and $n \in \NN$, we
set 
$\ord_{\varpi, \cX} (\omega) := \ord_{\varpi} (\tilde \omega) -n$.
Clearly this definition does not depend on the choice of
$\tilde \omega$.

\begin{lem}\label{ordre}Let $h : \cY \rightarrow \cX$ be a morphism between
flat formal $R$-schemes equidimensional
of relative dimension $d$. Let $\omega$ be in 
$\Omega^d_{\cX | R} (\cX)$ (resp. in $\Omega^d_{\cX_{K}} (\cX_K)$).
Let $y$ be a point in $\Gr (\cY) \setminus \Gr (\cY_{\rm sing})$
and assume $h (y)$ belongs to
$\Gr (\cX) \setminus \Gr (\cX_{\rm sing})$.
Then $$\ord_{\varpi} (h^{\ast} \omega) (y) =
\ord_{\varpi} (\omega) (h(y)) + \ord_{\varpi}(\Jac_h) (y),$$
resp.
$$\ord_{\varpi, \cY} (h^{\ast}_K \omega) (y) =
\ord_{\varpi, \cX} (\omega) (h(y)) + \ord_{\varpi}(\Jac_h) (y).$$
\end{lem}

\begin{proof}Follows directly from the definitions.
\end{proof}

\begin{def-theorem}Let $X$ be a smooth rigid variety over $K$, purely of 
dimension $d$. Let $\omega$ be 
a differential form in $\Omega^d_{X} (X)$.
\begin{enumerate}
\item[(1)]Let $\cX$ be a formal $R$-model
of $X$. Then the function $\ord_{\varpi, \cX} (\omega)$ is
exponentially integrable on
$\Gr (\cX)$ and the integral
$\int_{\Gr (\cX)} \LL^{- \ord_{\varpi, \cX} (\omega)} d\mu$ in $\widehat \kvar$
does not depend on the model $\cX$.
We denote it by
$\int_X \omega d\mu$.
\item[(2)]Assume furthermore $\omega$ is a gauge form, \textit{i.e.}
that
it generates $\Omega^d_{X}$ at every point of $X$,
and assume some open dense
formal subscheme  $\cU$ of $\cX$ is a
weak N{\'e}ron model of $X$.
Then the function $\ord_{\varpi, \cX} (\omega)$ takes
only a finite number of values and its fibres are stable cylinders.
Furthermore the integral
$\int_{\Gr (\cX)} \LL^{- \ord_{\varpi, \cX} (\omega)} d\tilde
\mu$ in $\kvar_{\rm loc}$
does not depend on the model $\cX$.
We denote it by
$\int_X \omega d \tilde \mu$.
\end{enumerate}
\end{def-theorem}

\begin{proof}Let us prove (2). Write $\omega = \varpi^{-n}
\tilde \omega$, with
$\tilde \omega$ in $\Omega^d_{\cX | R} (\cX)$ and $n \in \NN$.
Since $\cU$ is smooth, $\Omega^d_{\cU | R}$
is locally free of rank 1 
and $\tilde \omega \cO_U \otimes (\Omega^d_{\cU | R})^{-1}$
is isomorphic to a principal ideal sheaf $(f) \cO_U$, with
$f$ in $\cO_U$.
Furthermore, the function
$\ord_{\varpi, \cX} (\tilde \omega)$ coincides with
the 
function $\ord_{\varpi} (f)$
which to a point $\varphi$ of
$\Gr (\cU) = \Gr (\cX)$ associates
$\ord_{\varpi} (f (\varphi))$. The fibres of
$\ord_{\varpi} (f)$ are stable cylinders.
Since $\omega$ is a gauge form, $f$ induces an invertible function
on $X$, hence, by the maximum principle (cf. \cite{BGR}), the 
function $\ord_{\varpi} (f)$
takes only a finite number of values.
To prove 
that
$\int_{\Gr (\cX)} \LL^{- \ord_{\varpi, \cX} (\omega)} d\tilde
\mu$ in $\kvar_{\rm loc}$
does not depend on the model $\cX$, it is enough to consider the
case of another model $\cX'$ obtained from $\cX$
by
an admissible formal blow-up $h : \cX' \rightarrow \cX$. We may also assume
$\cX'$ contains as an 
open dense
formal subscheme 
a weak N{\'e}ron model $\cU'$ of $X$.
The equality
$$\int_{\Gr (\cX')} \LL^{- \ord_{\varpi, \cX'} (\omega)} d\tilde
\mu
=
\int_{\Gr (\cX)} \LL^{- \ord_{\varpi, \cX} (\omega)} d\tilde
\mu
$$
then follows from
Lemma \ref{ordre} and Theorem \ref{cvtilde}.
Statement (1) follows similarly from
Lemma \ref{ordre} and Theorem \ref{cv}.
\end{proof}

\begin{remark}
A situation where gauge forms naturally occur 
is that of reductive groups.
Let $G$ be a connected reductive group
over $k$.
B. Gross constructs in \cite{gross}, using Bruhat-Tits theory,
a differential 
form of top degree $\omega_G$ on $G$ which is defined up to
multiplication by a unit in $R$. One may easily check that
the differential form
$\omega_G$ induces a 
gauge form on the 
rigid $K$-group
$G^{\rm rig} := (G \widehat \otimes R)_K$.
\end{remark}

\begin{lem}\label{add}
Let $X$ be a smooth
rigid variety over $K$, purely of 
dimension $d$, and let $\omega$ be a gauge form on $X$.
Let $\cO= (O_i)_{i \in J}$ 
be a finite admissible covering and 
set $O_I := \cap_{i \in I} O_i$ for $I \subset J$.
Then
$$
\int_X \omega d \tilde \mu
=
\sum_{\emptyset \not= I \subset J}
(-1)^{|I| - 1} \int_{O_I} \omega_{| O_I} d \tilde \mu.
$$
If $\omega$ is only assumed to be a differential form in
$\Omega_X^d (X)$, then
$$
\int_X \omega d  \mu
=
\sum_{\emptyset \not= I \subset J}
(-1)^{|I| - 1} \int_{O_I} \omega_{| O_I} d \mu.
$$
\end{lem}

\begin{proof}Let us prove the first statement,
the proof of the second one being similar.
It is enough to consider the case $|J| = 2$.
Choose an $R$-model $\cX$ containing a
weak N{\'e}ron model $\cU$ of $X$
as an 
open dense
formal subscheme and such that the covering
$X = O_1 \cup O_2$ is induced from a covering
$\cX = \cO_1 \cup \cO_2$ by open formal subschemes.
It is sufficient to prove that
\begin{equation*}
\begin{split}
\int_{\Gr (\cX)} \LL^{- \ord_{\varpi, \cX} (\omega)} d\tilde \mu
= 
\int_{\Gr (\cO_1)} & \LL^{- \ord_{\varpi, \cO_1} (\omega_{| O_1})} d\tilde \mu
+
\int_{\Gr (\cO_2)} \LL^{- \ord_{\varpi, \cO_2} (\omega_{| O_2})} d\tilde \mu
\\
&-
\int_{\Gr (\cO_1 \cap \cO_2)} \LL^{- \ord_{\varpi, 
\cO_1 \cap \cO_2} (\omega_{| O_1 \cap O_2})} d\tilde \mu,\\
\end{split}
\end{equation*}
which follows from the fact that for every
open formal subscheme $\cO$ of $\cX$ the function
$\ord_{\varpi, \cX} (\omega)$ restricts to
$\ord_{\varpi, \cO} (\omega_{| \cO_K})$
on $\Gr (\cO)$ and the equalities 
$\Gr (\cX) = \Gr (\cO_1) \cup 
\Gr (\cO_2)$ and
$\Gr (\cO_1) \cap 
\Gr (\cO_2) = 
\Gr (\cO_1 \cap \cO_2)$ which follow from Proposition
\ref{glue}.
\end{proof}

\begin{prop}\label{prod}Let $X$ and $X'$ be smooth rigid $K$-varieties
purely of dimension $d$ and $d'$
and let 
$\omega$ and $\omega'$  be gauge forms on $X$ and $X'$.
Then
$$
\int_{X \times X'} \omega \times \omega' d \tilde \mu
=
\int_X \omega d \tilde \mu \times \int_{X'} \omega' d \tilde \mu.
$$
If $\omega$ and $\omega'$ are only assumed
to be differential forms in
$\Omega_X^d (X)$,
then
$$
\int_{X \times X'} \omega \times \omega' d  \mu
=
\int_X \omega d  \mu \times \int_{X'} \omega' d  \mu.
$$
\end{prop}

\begin{proof}Let us prove the first assertion,
the proof of the second one being similar.
Choose $R$-models $\cX$ and
$\cX'$ of $X$ and $X'$ respectively containing a
weak N{\'e}ron model $\cU$ of $X$ and $\cU'$ of $X'$
as an 
open dense
formal subscheme. Also write $\omega = \varpi^{-n} \tilde \omega$
and $\omega' = \varpi^{-n'} \tilde \omega'$, with 
$\tilde \omega$ and $\tilde \omega'$ in $\Omega^d_{\cX | R} (\cX)$
and 
$\Omega^{d'}_{\cX' | R} (\cX')$, respectively.
It enough to check that
$\tilde \mu (\ord_{\varpi, \cX \times \cX'} (\tilde \omega \times
\tilde \omega') = m)$
is equal to $\sum_{m'+ m'' = m} 
\tilde \mu (\ord_{\varpi, \cX} (\tilde \omega) = m') \times
\tilde \mu (\ord_{\varpi, \cX'} (\tilde \omega') = m'')$,
which follows the fact that on $\Gr (\cX \times \cX') \simeq
\Gr (\cX) \times \Gr (\cX') = \Gr (\cU) \times \Gr (\cU')$,
the functions
$\ord_{\varpi, \cX \times \cX'} (\tilde \omega \times
\tilde \omega')$ and
$\ord_{\varpi, \cX} (\tilde \omega) +
\ord_{\varpi, \cX'} (\tilde \omega')$ are equal.
\end{proof}

\subsection{Invariants for gauged smooth rigid varieties}
Let $d$ be an integer $\geq 0$.
We define $K_0 ({\rm GSRig}_K^d)$, the Grothendieck group of
gauged smooth rigid $K$-varieties of dimension $d$, as
follows: as an abelian group it
is the quotient of the free abelian group over symbols
$[X, \omega]$, with $X$ a
smooth rigid $K$-variety of dimension $d$
and $\omega$ a gauge form on $X$ by the  relations
$$
[X', \omega'] = [X, \omega]
$$
if there is an isomorphism $h : X' \rightarrow X$ with $h^{\ast} \omega 
= \omega'$, 
and 
$$
[X, \omega]
=
\sum_{\emptyset \not= I \subset J}
(-1)^{|I| - 1} [O_I, \omega_{| O_I}],
$$ when
$(O_i)_{i \in J}$
is a finite admissible covering of $X$, with the notation
$O_I := \cap_{i \in I} O_i$ for $I \subset J$.
One puts a graded ring structure
on $K_0 ({\rm GSRig}_K) := \oplus_d K_0 ({\rm GSRig}_K^d)$ by requiring
$$[X, \omega] \times 
[X', \omega'] := [X \times X', \omega \times \omega'].$$
Forgetting gauge forms, one
defines similarly 
$K_0 ({\rm SRig}_K^d)$,
the Grothendieck group of
smooth rigid $K$-varieties of dimension $d$, and the
graded ring
$K_0 ({\rm SRig}_K) := \oplus_d K_0 ({\rm SRig}_K^d)$.
There are natural forgetful morphisms
$$F : K_0 ({\rm GSRig}_K^d) \longrightarrow K_0 ({\rm SRig}_K^d)$$
and
$$F : K_0 ({\rm GSRig}_K) \longrightarrow K_0 ({\rm SRig}_K).$$

\begin{prop}\label{mainprop}
The assignment which to a gauged smooth rigid $K$-variety
$(X, \omega) $ associates 
$\int_X \omega d \tilde \mu$ 
factorizes uniquely as a ring morphism
$$
\tilde \mu : K_0 ({\rm GSRig}_K)
\rightarrow 
\kvar_{\rm loc}.
$$
\end{prop}

\begin{proof}This follows from Lemma \ref{add} and Proposition \ref{prod}.
\end{proof}

\subsection{A formula for $\int_X \omega d \tilde \mu$}
Let $X$ be a smooth rigid variety over $K$ of pure dimension $d$.
Let $\cU$ be
a weak N{\'e}ron model
of $X$ contained in some model $\cX$ of $X$
and let $\omega$ be a form in 
$\Omega^d_{\cX | R} (\cX)$ inducing a gauge form on $X$.
We denote by $U_0^i$, $i \in J$, the irreducible components
of the special fibre of $\cU$. By assumption,
each $U_0^{i}$ is smooth and $U_0^{i} \cap U_0^{j} = \emptyset$
for $i \not= j$.
We denote by $\ord_{U_0^{i}} (\omega)$ the unique integer $n$ such that
$\varpi^{-n} \omega$ generates $\Omega^d_{\cX | R}$ at the generic point of
$U_0^{i}$. More generally, if $\omega$ is a gauge form 
in $\Omega^d_{\cX_{K}} (\cX_K)$, 
we write $\omega = \varpi^{-n}
\tilde \omega$, with
$\tilde \omega$ in $\Omega^d_{\cX | R} (\cX)$ and $n \in \NN$, and we
set 
$\ord_{U_0^{i}} (\omega) := \ord_{U_0^{i}} (\tilde \omega) -n$.

\begin{prop}\label{form}Let $X$ be a smooth rigid variety
over $K$ of pure dimension $d$.
Let $\cU$ be
a weak N{\'e}ron model
of $X$ contained in some model $\cX$ of $X$
and let $\omega$ be a gauge
form in $\Omega^d_{\cX_{K}} (\cX_K)$.
With the above notations, 
we have
$$
\int_X \omega d \tilde \mu
=
\LL^{-d} \, \sum_{i \in J} \, [U_0^{i}] \, \LL^{- \ord_{U_0^{i}} (\omega)}
$$
in $\kvar_{\rm loc}$.
\end{prop}

\begin{proof}Denote by $\cU_0^{i}$ the irreducible component of $\cX$ 
with special fibre $U_0^{i}$.
Since $\Gr (\cX)$ is the disjoint union of the
sets $\Gr (\cU_0^{i})$, we may assume $\cX$ is a smooth irreducible
formal $R$-scheme of dimension $d$. Let
$\omega$ be a section of 
$\Omega^d_{\cX | R} (\cX)$ which generates
$\Omega^d_{\cX | R}$ at the generic point
of $\cX$ and induces a gauge form on the generic fibre.
Let us remark that the function
$\ord_{\varpi, \cX} (\omega)$ is identically equal to $1$
on $\Gr (\cX)$. Indeed, after shrinking $\cX$, we may write
$\omega = f \omega_0$ with $\omega_0$ a generator
$\Omega^d_{\cX | R}$ a every point and $f$ in $\cO_{\cX} (\cX)$.
By hypothesis 
$f$ is a unit at the generic point of $\cX$.
Assume at some point $x$ of $\Gr (\cX)$, 
$\ord_{\varpi} f (x) \geq 1$ ; it would
follow that the locus of $f = 0$ is non empty in $\cX$, which contradicts
the assumption that $\omega$
induces a gauge form on the generic fibre. Hence, we get
$\int_{\Gr (\cX)} \LL^{- \ord_{\varpi, \cX} (\omega)} \tilde d\mu = 
\LL^{-d} [X_0]$, and the result follows.
\end{proof}

\subsection{Application to Calabi-Yau varieties over $K$}
Let $X$ be a Calabi-Yau variety over $K$. By this we mean
a smooth projective algebraic variety over $K$ of pure dimension $d$
with $\Omega^d_X$ trivial.
We denote by $X^{\rm an}$
the rigid $K$-variety associated to
$X$. Since $X$ is proper, $X^{\rm an}$ is canonically isomorphic
to
the generic fibre
of the formal completion $X \widehat \otimes R$.
By GAGA (cf. \cite{lut} Theorem 2.8), $\Omega^d_{X^{\rm an}} (X^{\rm an})  = 
\Omega^d_{X} (X) = K$.

Now we can associate to any Calabi-Yau variety over $K$ a canonical
element in the ring
$\kvar_{\rm loc}$ which coincides with the class of 
the special fibre when $X$ has a model with good reduction.

\begin{theorem}\label{CY}Let $X$ be a Calabi-Yau variety over $K$,
let
$\cU$ be
a weak N{\'e}ron model
of $X^{\rm an}$  and let $\omega$ be 
a gauge form on $X^{\rm an}$.
We denote by $U_0^{i}$, $i \in J$, the irreducible components
of the special fibre of $\cU$ and set $\alpha (\omega) :=
\inf \ord_{U_0^{i}} (\omega)$.
Then the virtual variety
\begin{equation}\label{res}
[\overline X]
:=
\sum_{i \in J} \, [U_0^{i}] \, \LL^{\alpha (\omega)- \ord_{U_0^{i}} (\omega)}
\end{equation}
in $\kvar_{\rm loc}$
only depends on $X$. When $X$ has a proper smooth 
model with
good reduction
over $R$, $[\overline X]$ is equal to the class of the special fibre.
\end{theorem}

\begin{proof}Let $\omega$ be a
gauge form on $X^{an}$. By
Proposition \ref{form}, the right hand side of
(\ref{res}) is equal to 
$\LL^{d - \alpha (\omega)}\int_{X^{an}}
\omega d \tilde \mu$ which does not depend 
on $\omega$.
\end{proof}

In particular, we have the following analogue of Batyrev's result
on birational projective Calabi-Yau manifolds \cite{BCY}, \cite{arcs}.
\begin{cor}Let $X$ be a Calabi-Yau variety over $K$
and let $\cX$ and $\cX'$ be two proper and smooth
$R$-models of $X$ with special fibres
$\cX_0$ and $\cX'_0$.
Then $$[\cX_0] = [\cX'_0]$$ in
$\kvar_{\rm loc}$.\qed
\end{cor}

\begin{remark}
Calabi-Yau varieties over $k ((t))$, with $k$ of characteristic
zero have been considered in \cite{ks}.
\end{remark}

\subsection{A motivic Serre invariant for smooth rigid varieties}
We can now define the 
motivic Serre invariant for smooth rigid varieties.

\begin{theorem}\label{motivicserre}
There is a canonical ring morphism
$$\lambda : 
K_0 ({\rm SRig}_K) \longrightarrow \kvar_{\rm loc} / (\LL - 1)
\kvar_{\rm loc}$$
such that the 
diagram
\begin{equation*}\label{t}\xymatrix{
K_0 ({\rm GSRig}_K)
\ar[d]^{F} \ar[r]^<<<<<<<<<<<{\tilde \mu} & \kvar_{\rm loc} \ar[d]\\
K_0 ({\rm SRig}_K)
\ar[r]^<<<<<{\lambda}&\kvar_{\rm loc} / (\LL - 1)
\kvar_{\rm loc}
}
\end{equation*}
is commutative.
\end{theorem}

\begin{proof}Since any smooth rigid $K$-variety of dimension $d$
admits
a finite admissible covering by affinoids
$(O_i)_{i \in J}$, with $\Omega^d_{O_i}$ trivial, the morphism $F$ is
surjective.
Hence it is enough to show the following statement:
let $\cX$ be 
a smooth formal $R$-scheme
of relative dimension $d$
with $\Omega^d_{\cX | R}$ trivial, and let $\omega_1$ and $\omega_2$ be two
global sections of $\Omega^d_{\cX | R}$ inducing gauge forms 
on the generic fibre $\cX_K$,
then
$\int_{\Gr (\cX)} 
(\LL^{- \ord_{\varpi, \cX} (\omega_1)}-L^{- \ord_{\varpi, \cX} (\omega_2)})
d\tilde
\mu$
belongs to $(\LL - 1)
\kvar_{\rm loc}$.
To prove this,  we take $\omega_0$ a
global section of $\Omega^d_{\cX | R}$ such that
$\Omega^d_{\cX | R} \simeq \omega_0 \cO_{\cX}$.
If $\omega$ is any global section of $\Omega^d_{\cX | R}$, write $\omega = f 
\omega_0$ with $f$ in  $\cO_{\cX} (\cX)$.
By 
the maximum principle, the 
function $\ord_{\varpi} (f)$
takes only a finite number of values on $\Gr (\cX)$.
It follows we may write $\Gr (\cX)$ as a disjoint union of the subsets
$\Gr (\cX)_{\ord_{\varpi} (f) = n}$ where $\ord_{\varpi} (f)$ takes 
the value $n$. These subsets
are stable cylinders and only a finite number of them are non empty.
Hence the equality
$$
\int_{\Gr (\cX)} 
(\LL^{- \ord_{\varpi, \cX} (\omega)}-L^{- \ord_{\varpi, \cX} (\omega_0)})
d\tilde
\mu
=
\sum_n (\LL^{-n} - 1) \tilde \mu  (\Gr (\cX)_{\ord_{\varpi} (f) = n})
$$
holds in $\kvar_{\rm loc}$ and the statement follows.
\end{proof}

\begin{remark}The ring $\kvar_{\rm loc} / (\LL - 1)
\kvar_{\rm loc}$ is much smaller than
$\kvar_{\rm loc}$, but still quite
large. Let $\ell$ be a prime number 
distinct from the characteristic of $k$. Then the {\'e}tale
$\ell$-adic
Euler characteristic with compact supports
$ X \mapsto \chi_{c, \ell}  (X) := \sum (-1)^i \dim H^i_{c, \text{\'et}}
(X, \QQ_{\ell})$
induces a ring morphism $\chi_{c, \ell} : \kvar_{\rm loc} / (\LL - 1)
\kvar_{\rm loc} \rightarrow \ZZ$.
Similarly, assume there is a natural morphism
$H : \kvar_{\rm loc} \rightarrow \ZZ [u, v]$
which to the class of
a variety $X$ over $k$ assigns its Hodge polynomial $H (X)$ for
de Rham cohomology with compact support. Such a morphism is known
to exist when $k$ is of characteristic zero.
Then if one sets $H_{1/2} (X) (u) := H (X) (u, u^{-1})$
one gets a morphism $H_{1/2} : \kvar_{\rm loc} / (\LL - 1)
\kvar_{\rm loc} \rightarrow \ZZ [u]$, since $H (\AA^1_k) = uv$.
\end{remark}

\begin{theorem}\label{invner}
Let $X$ be a smooth rigid variety over $K$ of pure dimension $d$.
Let $\cU$ be
a weak N{\'e}ron model
of $X$ and denote by
$U_0$ its special fibre.
Then
$$
\lambda ([X]) = [U_0]
$$
in $\kvar_{\rm loc} /(\LL - 1)\kvar_{\rm loc}$.

In particular,  the class of
$[U_0]$ in $\kvar_{\rm loc} /(\LL - 1)\kvar_{\rm loc}$
does not depend on  the weak N{\'e}ron model $\cU$.
\end{theorem}

\begin{proof}By taking an appropriate admissible cover, we may assume
there exists a gauge form on $X$, in which case the result follows from
Proposition \ref{form},
since
$[U_0] = \sum_{i \in J} [U_0^i]$.
(In fact, one can also prove
Theorem \ref{motivicserre} that way, but we prefered to give a proof
which is quite parallel to that of Serre in \cite{serre}.)
\end{proof}

\subsection{Relation with $p$-adic integrals on
compact locally analytic varieties}
Let $K$ be a local field with residue field $k = \FF_q$.
Let us consider the Grothendieck group 
$K_0 ({\rm SLocAn}^d_K)$ of compact locally analytic
smooth varieties over $K$ of pure dimension $d$, which is defined similarly as
$K_0 ({\rm SRig}^d_K)$ replacing smooth rigid varieties by
compact locally analytic
smooth varieties and finite admissible covers by finite covers.
Also a nowhere vanishing
locally analytic $d$-form on a smooth 
compact locally analytic variety $X$ of pure dimension $d$ will
be called a gauge form on $X$, and one defines
the Grothendieck group 
$K_0 ({\rm GSLocAn}^d_K)$  of gauged
compact locally analytic
smooth varieties over $K$ of pure dimension $d$ similarly as
$K_0 ({\rm GSRig}^d_K)$. There are canonical
forgetful morphisms
$F : K_0 ({\rm SRig}^d_K) \rightarrow K_0 ({\rm SLocAn}^d_K)$
and
$F : K_0 ({\rm GSRig}^d_K) \rightarrow K_0 ({\rm GSLocAn}^d_K)$
induced from the functor which to 
a rigid variety (resp. gauged variety) associates the underlying
locally analytic variety (resp. gauged variety).
If $(X, \omega)$ 
is a  gauged
compact locally analytic
smooth variety, the $p$-adic integral $\int_X |\omega|$ belongs to
$\ZZ [q^{-1}]$ (cf. \cite{serre}), and by additivity of $p$-adic integrals
one gets a morphism
${\rm int}_p : K_0 ({\rm GSLocAn}^d_K) \rightarrow \ZZ [q^{-1}]$.

On the other hand, there is a canonical morphism
$N : \kvar \rightarrow \ZZ$ which to the class of a $k$-variety $S$
assigns the number of points of $S (\FF_q)$, and which
induces a morphism
$N : \kvar_{\rm loc} \rightarrow \ZZ [q^{-1}]$.
We shall also denote by $N$ the induced morphism
$\kvar_{\rm loc} /(\LL - 1)\kvar_{\rm loc}  \rightarrow 
\ZZ [q^{-1}] / (q - 1)\ZZ [q^{-1}] \simeq  \ZZ  / (q - 1)\ZZ$.

\begin{prop}\label{compint}Let $K$ be a local field with residue field $k = \FF_q$.
Then the diagram
\begin{equation*}
\label{ttt}\xymatrix{
K_0 ({\rm GSRig}^d_K)
\ar[d]^{F} \ar[r]^<<<<{\tilde \mu} & \kvar_{\rm loc}  \ar[d]^{N}\\
K_0 ({\rm GSLocAn}^d_K)
\ar[r]^<<<<<{{\rm int}_p}&\ZZ [q^{-1}]
}
\end{equation*}
is commutative.
\end{prop}

\begin{proof}One reduces to showing the following: 
let $\cX$ be a smooth formal $R$-scheme of dimension $d$ and let $f$
be a function in $\cO_{\cX} (\cX)$ which induces a non vanishing function
on $\cX_K$, then $$
N (\int_{\Gr(\cX)} \LL^{- \ord_{\varpi} (f)} d\tilde \mu)
=
\int_{\cX (R)} q^{- \ord_{\varpi} (f)} d\tilde \mu_p,
$$
with $d\tilde \mu_p$ the $p$-adic measure on $\cX (R)$.
It is enough to check that
$N (\tilde \mu (\ord_{\varpi} (f) = n))$ is equal to the $p$-adic measure 
of the set of points $x$ of $\cX (R)$ with $\ord_{\varpi} (f) (x) = n$,
which follows from Lemma \ref{compmes}.
\end{proof}

\begin{lem}\label{compmes}Let $K$ be a local field with residue field $k = \FF_q$. Let $\cX$ be a smooth formal $R$-scheme of dimension $d$.
Let $A$ be a (stable) cylinder 
in $\Gr (\cX)$.
Then $N (\tilde \mu (A))$ is equal to the $p$-adic volume of
$A \cap \Gr (\cX) (k)$.
\end{lem}

\begin{proof}Write $A = \pi_n^{-1}(C)$, with $C$ a constructible
subset if $\Gr_n (\cX)$. By definition 
$\tilde \mu (A) = \LL^{ -d (n + 1)} [C]$. On the other hand 
$\cX$ being smooth, the morphism $A \cap \Gr (\cX) (k) \rightarrow C (k)$
is surjective and its fibres are balls of radius $q^{-d (n + 1)}$.
It follows that the $p$-adic volume of $A \cap \Gr (\cX) (k)$
is equal to $|C (k)| q^{-d (n + 1)}$.
\end{proof}

Let us now explain the relation with the work of Serre
in \cite{serre}. 
Serre shows in \cite{serre}
that any compact locally analytic
smooth variety over $K$ of pure dimension $d$
is isomorphic to $r B^d$, with $r$
an integer $\geq 1$ and $B^d$ the unit ball
of dimension $d$ and that, futhermore,
$r B^d$ is isomorphic to $r' B^d$ if and only if $r$ and $r'$ are congruent
modulo $q - 1$. We shall denote by $s (X)$ the class of $r$ in
$\ZZ /(q - 1) \ZZ$. 
It follows from
Serre's results that $s$ induces an isomorphism
$s : K_0 ({\rm SLocAn}^d_K) \rightarrow \ZZ /(q - 1) \ZZ$ and
that the diagram
\begin{equation*}
\xymatrix{
K_0 ({\rm GSLocAn}^d_K)
\ar[d] \ar[r]^<<<<{{\rm int}_p} & \ZZ [q^{-1}] \ar[d]\\
K_0 ({\rm SLocAn}^d_K)
\ar[r]^<<<<<{s}&\ZZ /(q - 1) \ZZ
}
\end{equation*}
is commutative. The following result then follows from Proposition
\ref{compint}.

\begin{cor}Let $K$ be a local field with residue field $k = \FF_q$.
Then the diagram
\begin{equation*}
\label{tt}\xymatrix{
K_0 ({\rm SRig}^d_K)
\ar[d]^{F} \ar[r]^<<<<{\lambda} & \kvar_{\rm loc} /
(\LL - 1)\kvar_{\rm loc} \ar[d]^N\\
K_0 ({\rm SLocAn}^d_K)
\ar[r]^<<<<<<<<<<<{s}&\ZZ /(q - 1) \ZZ
}
\end{equation*}
is commutative. \qed
\end{cor}

\section{Essential components of weak N{\'e}ron models}\label{sec5}
\subsection{Essential components and the Nash problem}
Since we shall proceed 
by analogy with \cite{nash}, let us begin by recalling some material from
that paper.
We assume in this subsection that 
$k$ is of characteristic zero and that $R = k [[\varpi]]$.
For $X$ an algebraic variety over $k$, we denote by $\cL (X)$ its arc space
as defined in \cite{arcs}. In fact, in the present \S\kern .15em,
we shall use 
notations and results from \cite{arcs}, even when they happen to be
 special cases
of ones in this paper.
As remarked in  \ref{arcspace},
$\cL (X) = \Gr (X \widehat \otimes R)$
and there are natural morphisms
$\pi_n : \cL (X) \rightarrow \cL_n (X)$
with 
$\cL_n (X) = \Gr_n (X \otimes R_n)$.

By a desingularization
of a variety $X$ we mean a proper and birational morphism
$$
h : Y \longrightarrow X,
$$
with $Y$ a smooth variety, inducing an isomorphism
between $h^{-1} (X \setminus X_{\rm sing}) $
and $X \setminus X_{\rm sing} $ (some authors omit
the last condition).

Let $h : Y \rightarrow X$ be a desingularization of $X$ and let
$D$ be an irreducible component of
$h^{-1} (X_{\rm sing}) $
of codimension 1 in $Y$.
If $h' : Y' \rightarrow X$ is another 
desingularization of $X$, the birational map
$\phi: h'{}^{-1} \circ h : Y \dashrightarrow Y'$
is defined at the generic point $\xi$
of $D$, since $h'$ is proper, hence we can define the image of 
$D$ in $Y'$ as the closure of $\phi (\xi)$ in $Y'$.
One says  that $D$ is an essential divisor
with respect to $X$, if, for every desingularization
$h' : Y' \rightarrow X$ of $X$ the image of $D$ in $Y'$ is a divisor
and that $D$ is an essential component 
with respect to $X$, if, for every desingularization
$h' : Y' \rightarrow X$ of $X$ the image of $D$ in $Y'$
is an irreducible component of 
$h'{}^{-1} (X_{\rm sing})$.
In general, if $D$ is an irreducible component of
$h^{-1} (X_{\rm sing}) $, we say 
$D$ is an essential component with respect to $X$, if
there exists
a proper birational morphism
$p : Y'\rightarrow Y$, with $Y'$ smooth, and a
divisor $D'$ in $Y'$ such that $D'$ is an essential component 
with respect to $X$ and $p(D') = D$.
It follows from
the definitions and Hironaka's Theorem
that
essential components of different resolutions of the same variety $X$
are in natural bijection, hence we may denote by
$\tau (X)$ the number of  essential components in any resolution of $X$.

Let $W$ be a constructible subset of an algebraic variety $Z$.
We say that $W$ is irreducible in $Z$
if the Zariski closure $\overline W$
of $W$ in $Z$ is irreducible.

In general let $\overline W = \cup_{1 \leq i \leq n} W'_{i}$ 
be the decomposition of $\overline W$ into irreducible components.
Clearly  $W_{i} := W'_{i} \cap W$ is non empty, irreducible in $Z$,
and its closure in $Z$ is equal
to $W'_{i}$. We call the $W_{i}$'s the irreducible components of
$W$ in $Z$.

Let $E$ be a locally closed subset of $h^{-1} (X_{\rm sing})$.
We denote by $Z_{E}$ the set of arcs in $\cL (Y)$ whose origin
lies on $E$ but which are not contained in $E$. In other words
$Z_{E} = \pi^{-1}_{0} (E) \setminus \cL (E)$.
Let us remark that if $E$ is smooth and connected, $\pi_n (Z_{E})$
is constructible and irreducible in 
$\cL_n (Y)$.
Now we set $N_{E} := h (Z_{E})$. Since $\pi_n (N_{E})$
is the image of $\pi_n (Z_{E})$ under the morphism
$\cL_n (Y) \rightarrow \cL_n (X)$ induced by $h$, it follows that
$\pi_n (N_{E})$
is constructible and irreducible in 
$\cL_n (Y)$.

The following result, proved in \cite{nash}, follows easily from the above 
remarks and Hironaka's
resolution of singularities:

\begin{prop}[Nash \cite{nash}]\label{nashprop}Let $X$ be an algebraic variety over $k$, a field of
characteristic zero.
Set 
$\cN (X) := \pi^{-1}_{0} (X_{\rm sing}) \setminus \cL (X_{\rm sing})$.
For every $n \geq 0$, $\pi_n (\cN (X))$ is a constructible subset
of $\cL_n (X)$.
Denote by
$W^1_n$, \dots $W^{r (n)}_n$  the irreducible components of
$\pi_n (\cN (X))$. The mapping $n \mapsto r(n)$ is nondecreasing and bounded
by the number $\tau (X)$ of essential components occuring
in  a resolution of $X$.
\end{prop}

Up to  renumbering, we may assume that
$W^{i}_{n +1}$ maps to 
$\overline {W^{i}_{n}}$ for $n \gg 0$.
Let us call the family
$(W^{i}_{n})_{n \gg 0}$ a Nash family.
Nash shows furthermore
that for every Nash family $(W^{i}_{n})_{n \gg0}$
there exists a unique
essential component $E$ in  a given resolution $h : Y \rightarrow
X$ of $X$ such that 
$\overline {\pi_n (N_{E})} = \overline {W^{i}_{n}}$ for $n \gg 0$.

Now, we can formulate 
the Nash problem:

\begin{problem}[Nash \cite{nash} p.36]\label{nashprob}
Is there always a corresponding Nash family
for an essential component? In general, how completely do the essential
components correspond to Nash families? 
What is the relation between $\tau (X)$
and $\lim r (n)$?
\end{problem}

Let $W$ be a constructible subset
of some variety $X$. We denote  the
supremum of the dimension of the irreducible components
of the closure of $W$ in $Y$ by ${\rm dim} \, W$.

Let $h : Y \rightarrow X$ be a proper birational morphism
with $Y$ a smooth variety. Let $E$ be a codimension 1 irreducible
component of the exceptional locus   of $h$ in $Y$.
We denote by $\nu (E) - 1$ the length of
$\Omega^d_{Y} / h^{\ast} \Omega^d_{X}$ at the generic point
of $E$. Here $\Omega^d_{X}$ denotes the
$d$-th exterior power of the sheaf $\Omega^1_{X}$ of differentials on $X$.

\begin{prop}\label{codim2}Let $X$ be a variety of pure dimension $d$
over $k$, a field of characteristic 0.
Let $h : Y \rightarrow X$ be a proper birational morphism
with $Y$ a smooth variety and let
$U$ be a non empty
open subset
of
a codimension 1 irreducible
component $E$ of the exceptional locus   of $h$ in $Y$.
Then 
$${\rm dim} \, \pi_{n} (N_U)
=
(n + 1) \, d -  \nu (E)
$$
for $n \gg 0$.
\end{prop}

\begin{proof}By Theorem 6.1 of \cite{arcs},
the image of
$[\pi_{n} (N_U)] \LL^{- (n + 1) d}$ in $\widehat{\kvar}$
converges to  $\mu (N_U)$ in 
$\widehat{\kvar}$. 
Since ${\rm dim} \, \pi_{n} (N_U) \leq 
(n+ 1)d$
by Lemma 4.3 of \cite{arcs}, one deduces the fact that 
${\rm dim} \, \pi_{n} (N_U)
-
(n + 1) \, d$ has a limit. To conclude we 
first remark that
$\overline{\pi_{n} (N_{U})} =
\overline{\pi_{n} (N_{E})}$
for any non empty open subset
$U$ in $E$. Hence we may
assume that $(h^{\ast} \Omega^d_{X}) / {\rm torsion}$ is
locally free
on a neighborhood of  $U$. It then follows from
Proposition 6.3.2 in \cite{arcs}, or rather from its proof,
that
$$
\mu (N_{U}) = \LL^{-d} [U] (\LL - 1)
\sum_{\ell \geq 1} \LL^{- \ell \nu (E)}
$$
in $\widehat{\kvar}$. Hence $\mu (N_{U})$ belongs
to
$F^{\nu (E)}$ and not to $F^{\nu (E) + 1}$,
and the result follows.
\end{proof}

\subsection{Essential components of weak N{\'e}ron models}We shall return now
to the setting of the present paper.
We shall fix a flat formal $R$-scheme $\cX$
of relative dimension $d$ with smooth
generic fibre $\cX_K$. By a weak N{\'e}ron model of $\cX$, we shall mean
a weak N{\'e}ron model $\cU$ of 
$\cX_K$ together with a morphism $h : \cU \rightarrow \cX$ inducing
the inclusion $\cU_K \hookrightarrow \cX_K$.
As before we shall  denote by $U_0^{i}$, $i \in J$, the irreducible components
of the special fibre of $\cU$.
Let $\xi^i$ denote the generic point of $U_0^{i}$.
We shall say $U_0^{i}$ is an essential component with respect to $\cX$
if, for every 
weak N{\'e}ron model $\cU'$ of 
$\cX$, the Zariski closure of
${\pi_{0, \cU'} (\pi_{0, \cU}^{-1} (\xi^i))}$
is an irreducible component of the special fibre of $\cU'$.
Note that being an essential component is a property relative to $\cX$.
By their very definition, essential components in 
different weak N{\'e}ron models of 
$\cX$ are in natural bijection.

We have the following analogue of Proposition \ref{nashprop}.

\begin{prop}\label{annashprop}Let $\cX$ be a flat formal $R$-scheme 
of relative dimension $d$ with smooth
generic fibre $\cX_K$.
Denote by
$W^1_n$, \dots $W^{r (n)}_n$  the irreducible components of the constructible
subset 
$\pi_n (\Gr (\cX))$ of $\Gr_n (\cX)$.
The mapping $n \mapsto r(n)$ is nondecreasing and bounded
by the number $\tau (\cX)$ of essential components occuring in a
weak N{\'e}ron model of $\cX$.
\end{prop}

\begin{proof}Clearly
the mapping $n \mapsto r(n)$ is nondecreasing.
Let $h : \cU \rightarrow \cX$ be a weak N{\'e}ron model of $\cX$, with
irreducible components $\cU^i$, $i \in J$.
Since $\cU^i$ is smooth and irreducible, $\pi_{n, \cU} (\Gr (\cU^i))$
is also smooth and irreducible, hence the Zariski closure
of
$h (\pi_{n, \cU} (\Gr (\cU^i))) = \pi_{n, \cX} (\Gr (\cU^i))$
in $\Gr_n (\cX)$
is irreducible. Since $\Gr (\cX)$ is the union of the subschemes
$\Gr (\cU^i)$, it follows that $r (n)$ is bounded
by $|J|$.
Now if $\cU^i$ is not an essential component, there exists some
weak N{\'e}ron model of $\cX$, $h' : \cU' \rightarrow \cX$,
such that, if we denote by $W^i$
the image  of 
$\Gr (\cU^i)$ in $\Gr (\cU')$,
$\pi_{n, \cU'} (W^i)$ is contained in the Zariski closure
of
$\pi_{n, \cU'} (\Gr (\cU') \setminus W^i)$.
It follows that 
$\pi_{n, \cX} (\Gr (\cU^i))$ is contained in the closure of
$\pi_{n, \cX} (\Gr (\cU) \setminus \Gr (\cU^i))$.
The bound $r (n) \leq \tau (\cX)$ follows.
\end{proof}

As previously,
we may, up to  renumbering, assume that
$W^{i}_{n +1}$ maps to 
$\overline {W^{i}_{n}}$ for $n \gg 0$.
We shall still call the family
$(W^{i}_{n})_{n \gg 0}$ a Nash family.
Let $\xi^i_n$ be the generic point of
$\overline {W^{i}_{n}}$. By construction $\xi^i_{n +1}$ 
maps to $\xi^i_{n}$ 
under the truncation morphism
$\Gr_{n +1} (\cX) \rightarrow \Gr_{n} (\cX)$,
hence to the
inverse system $(\xi^{i}_{n})_{n \gg 0}$ corresponds a point
$\xi^{i}$ of $\Gr (\cX)$.
Let $h : \cU \rightarrow \cX$ be 
weak N{\'e}ron model of $\cX$ with irreducible components
$\cU^j$, $j \in J$. 
There is a unique irreducible component
$\cU^{j (i)}$ such that 
the point $\xi^i$ belongs to $\Gr (\cU^{j (i)})$.
Furthermore, 
$\overline {h (\pi_n (\Gr (\cU^{j (i)}))} =
\overline {W^{i}_{n}}$ for $n \gg 0$ and it follows from
the proof of 
Proposition \ref{annashprop}
that 
$U^{j (i)}_0$ is essential.

We have the following analogue of 
the Nash problem:

\begin{problem}\label{annashprob}
Is there always a corresponding Nash family
for an essential component? In general, how completely do the essential
components correspond to Nash families?
What is the relation between $\tau (\cX)$
and $\lim r (n)$?
\end{problem}

For $E$ a locally closed subset of the special fibre of $\cU$,
we set $Z_E := \pi_{0, \cU}^{-1} (E))$ and $N_E := h (Z_E)$.
Remark that, for every $n$, $\pi_n (Z_E)$ and
$\pi_n (N_E)$ are constructible subsets of $\Gr_n (\cU)$ and
$\Gr_n (\cX)$, respectively.
Indeed, $\pi_n (Z_E)$ is constructible since $\cU$ is smooth, hence
$\pi_n (N_E) = h (\pi_n (Z_E))$ also.
We denote by $\nu (U_0^{i}) - 1$ the length of
$\Omega^d_{\cU | R} / h^{\ast} \Omega^d_{\cX | R}$ at the generic point
$\xi^i$
of $U_0^{i}$.

We also have the following analogue of
Proposition \ref{codim2}.

\begin{prop}\label{codim3}Let $\cX$ be a flat formal $R$-scheme 
of relative dimension $d$ with smooth
generic fibre $\cX_K$.
Let $h : \cU \rightarrow \cX$ be a weak N{\'e}ron model of $\cX$
and let $E$ be an open dense subset of an irreducible component $U_0^{i}$
of the special fibre of $\cU$.
Then 
$${\rm dim} \, \pi_{n} (N_E)
=
(n + 1) \, d -  \nu (U_0^{i})
$$
for $n \gg 0$.
\end{prop}

\begin{proof}Fix an integer $e \geq 0$.
By Lemma \ref{cv2}, for $n \gg e$,
\begin{equation*}
\begin{split}\dim \pi_n (N_E \cap \Gr^{(e)} \cX)
&= \dim \pi_n (h^{-1} (N_E \cap \Gr^{(e)} \cX)) -  \nu (U_0^{i})\\
&= (n +1) d -  \nu (U_0^{i}).\\
\end{split}
\end{equation*}
On the other hand, it follows from Lemma \ref{dim2}
that
$$
\dim \pi_n (N_E \cap (\cX \setminus \Gr^{(e)} \cX))
<
(n +1) d -  \nu (U_0^{i}),
$$
when $n \gg e \gg \nu (U_0^{i})$.
\end{proof}

\bibliographystyle{amsplain}

\end{document}